%	fichier de macros mathematiques 
%
\parindent=0mm
\font \fiverm=cmr5
\font \sixrm=cmr6
\font \sevenrm=cmr7
\font \eightrm=cmr8
\font \ninerm=cmr9

\font \ninebf=cmbx9
\font \bigbf=cmbx10 scaled \magstep1
\font \Bigbf=cmbx10 scaled \magstep2

%		GOTHIQUE
\font \tengoth=eufm10
\font \sevengoth=eufm7
\font \fivegoth=eufm5

\newfam\gothfam
\textfont \gothfam=\tengoth
\scriptfont \gothfam=\sevengoth
\scriptscriptfont \gothfam=\fivegoth

%
%		PETIT ROMAIN
%
\newfam\srmfam
\textfont \srmfam=\eightrm
\scriptfont \srmfam=\sixrm
\scriptscriptfont \srmfam=\fiverm

%		GOTHIQUE GRAS
\font \tengothb=eufb10
\font \sevengothb=eufb7
\font \fivegothb=eufb5

\newfam\gothbfam
\textfont \gothbfam=\tengothb
\scriptfont \gothbfam=\sevengothb
\scriptscriptfont \gothbfam=\fivegothb

%		CARACTERES MATHEMATIQUES
\font \tenmath=msbm10
\font \sevenmath=msbm7
\font \fivemath=msbm5

\newfam\mathfam
\textfont \mathfam=\tenmath
\scriptfont \mathfam=\sevenmath
\scriptscriptfont \mathfam=\fivemath
\def\math{\fam\mathfam\tenmath}
%
%		PAGINATION
%
\def\titre#1{\centerline{\Bigbf #1}\nobreak\nobreak\vglue 10mm\nobreak}

\def\paragraphe#1{\bigskip\goodbreak {\bigbf #1}\nobreak\vglue 12pt\nobreak}
\def\alinea#1{\medskip\allowbreak{\bf#1}\nobreak\vglue 9pt\nobreak}
\def\ssq{\smallskip\qquad}

%
%		THEOREMES ET PROPOSITIONS
%
\def\th#1{\bigskip\goodbreak {\bf Theor\`eme #1.} \par\nobreak \sl }
\def\prop#1{\bigskip\goodbreak {\bf Proposition #1.} \par\nobreak \sl }
\def\lemme#1{\bigskip\goodbreak {\bf Lemme #1.} \par\nobreak \sl }
\def\cor#1{\bigskip\goodbreak {\bf Corollaire #1.} \par\nobreak \sl }
\def\dem{\bigskip\goodbreak \it D\'emonstration. \rm}
\def\ndem{\bigskip\goodbreak \rm}
\def\qed{\par\nobreak\hfill $\bullet$ \par\goodbreak}
%
%		DIVERS
%
\def\uple#1#2{#1_1,\ldots ,{#1}_{#2}}
\def\corde#1#2-#3{{#1}_{#2},\ldots ,{#1}_{#3}}

\def\ordcorde#1#2-#3{{#1}_{#2} \le \cdots \le {#1}_{#3}}
\def\strictordcorde#1#2-#3{{#1}_{#2} < \cdots < {#1}_{#3}}
\def \restr#1{\mathstrut_{\bigl |}\raise-8pt\hbox{$\scriptstyle #1$}}
\def \srestr#1{\mathstrut_{\textstyle |}\raise-6pt\hbox{$\scriptscriptstyle #1$}}
\def \inver{^{-1}}
\def\dbar{d\!\!\hbox to 4.5pt{\hfill\vrule height 5.5pt depth -5.3pt
	width 3.5pt}}

\def\frac#1#2{{\textstyle {#1\over #2}}}

\def\R{{\math R}}

\def\N{{\math N}}

\def\fleche#1{\mathop{\hbox to #1 mm{\rightarrowfill}}\limits}
\def\gfleche#1{\mathop{\hbox to #1 mm{\leftarrowfill}}\limits}
\def\inj#1{\mathop{\hbox to #1 mm{$\lhook\joinrel$\rightarrowfill}}\limits}
\def\ginj#1{\mathop{\hbox to #1 mm{\leftarrowfill$\joinrel\rhook$}}\limits}
\def\surj#1{\mathop{\hbox to #1 mm{\rightarrowfill\hskip 2pt\llap{$\rightarrow$}}}\limits}
\def\gsurj#1{\mathop{\hbox to #1 mm{\rlap{$\leftarrow$}\hskip 2pt \leftarrowfill}}\limits}
%
%		FONTES
%
\def \g#1{\hbox{\tengoth #1}}

\def\Cal #1{{\cal #1}}
%
%		OPERATEURS MATHEMATIQUES
%

\def \mop#1{\mathop{\hbox{\rm #1}}\nolimits}
\def \smop#1{\mathop{\hbox{\sevenrm #1}}\nolimits}

\def \mopl#1{\mathop{\hbox{\rm #1}}\limits}

%
%		BIBLIOGRAPHIE
%
\def \bib #1{\null\medskip \strut\llap{[#1]\quad}}

\magnification=\magstep1
\parindent=0cm
\def\titre#1{\centerline{\Bigbf #1}\vskip 16pt}
\def\paragraphe#1{\bigskip {\bigbf #1}\vskip 12pt}
\def\alinea#1{\medskip{\bf #1}\vskip 6pt}
\def\ssq{\smallskip\qquad}

\let\wt=\widetilde

\titre {Distributions \`a support compact}
\vskip -4mm
\titre {et repr\'esentations unitaires}
\vskip 8mm
\centerline{Dominique MANCHON}
\centerline{Institut Elie Cartan - UMR 9973 du CNRS}
\centerline{BP 239, 54506 Vandoeuvre cedex}
\centerline{Adresse \'electronique : {\tt manchon@iecn.u-nancy.fr}}
\vskip 12mm
{\ninebf R\'esum\'e}~: {\ninerm Dans cet article nous pr\'ecisons deux notions introduites par Roger Howe [Hw] pour \'etudier les repr\'esentations unitaires des groupes de Lie~: les repr\'esentations unitaires fortement tra\c cables et le front d'onde d'une repr\'esentation unitaire. Nous montrons que pour toute distribution $\varphi$ \`a support compact sur un groupe de Lie connexe dont le front d'onde ne rencontre pas l'oppos\'e du front d'onde de la repr\'esentation $\pi$ l'op\'erateur $\pi(\varphi)$ est r\'egularisant, g\'en\'eralisant ainsi un r\'esultat de [M2]. De plus sous les m\^emes hypoth\`eses $\pi(\varphi)$ est \`a trace si $\pi$ est fortement tra\c cable. Dans le cas o\`u la repr\'esentation $\pi$ est irr\'eductible et associ\'ee par la m\'ethode des orbites \`a une orbite $\Omega$ ferm\'ee et temp\'er\'ee, nous montrons qu'elle est fortement tra\c cable et nous \'etendons la formule des caract\`eres aux op\'erateurs $\pi(\varphi)$ pour les distributions $\varphi$ \`a support compact dont le front d'onde v\'erifie la condition de transversalit\'e ci-dessus.}
\bigskip
{\ninebf Abstract}~: {\ninerm We develop in this article two notions due to R. Howe [Hw] related to unitary representations of Lie groups, namely strong trace class representations as well as the wave front set of a unitary representation. Generalizing a previous theorem [M2] we show that for any compactly supported distribution $\varphi$ on a Lie connected group the wave front set of which does not intersect the opposite of the wave front set of the representation $\pi$ the operator $\pi(\varphi)$ is a smoothing operator, moreover a trace class operator if $\pi$ is of strong trace class. In case the representation $\pi$ is irreducible and associated to a closed and tempered coadjoint orbit $\Omega$ we show that it is of strong trace class, and we extend the character formula to the operators $\pi(\varphi)$ for any compactly supported distribution $\varphi$ satisfying the above transversality condition.}
\paragraphe{Introduction}
\qquad
Dans la correspondance de Kirillov entre orbites et repr\'esentations le caract\`ere de la repr\'esentation $\pi$ est donn\'e par une int\'egrale sur l'orbite coadjointe $\Omega$ correspondante. Nous nous proposons d'\'etendre cette formule des caract\`eres \`a certains op\'erateurs $\pi(\varphi)$ o\`u $\varphi$ est une distribution \`a support compact sur le groupe.
\ssq
Le probl\`eme peut se formuler bri\`evement comme suit~: soit $G$ un groupe de Lie connexe d'alg\`ebre de Lie $\g g$, et soit $\pi$ une repr\'esentation unitaire irr\'eductible de $G$ associ\'ee \`a une orbite coadjointe $\Omega$ ferm\'ee, admettant un caract\`ere-distribution donn\'e par une int\'egrale orbitale ``\`a la Kirillov''~:
$$\mop{Tr}\pi(\varphi)=\kappa(\pi)\int_\Omega \Cal F\bigl(P_\Omega\inver .j_G.(\varphi\circ\exp)\bigr)(\omega)\, d\beta_\Omega(\omega)$$
o\`u $\kappa(\pi)$ est un entier, $j_G$ le jacobien de l'exponentielle et $P_\Omega$ est une fonction analytique. Si $\varphi$ est une distribution \`a support compact sur $G$ on peut d\'efinir l'op\'erateur $\pi(\varphi)$ sur les vecteurs $C^\infty$ de la repr\'esentation. Or il r\'esulte des travaux de J-Y. Charbonnel [Cha 2,3] que l'orbite $\Omega$, \'etant ferm\'ee, est temp\'er\'ee. Donc si le front d'onde $WF(\varphi)$ en l'\'el\'ement neutre ne rencontre pas le c\^one asymptote \`a $\Omega$ et si le support de $\varphi$ est assez petit, l'int\'egrale orbitale converge. Peut-on pour autant en d\'eduire que l'op\'erateur $\pi(\varphi)$ est \`a trace et \'etendre la formule des caract\`eres aux op\'erateurs de ce type?  
\ssq
Nous donnons une r\'eponse affirmative \`a cette question en repla\c cant ce probl\`eme dans un cadre g\'en\'eral trac\'e par Roger Howe [Hw] autour des notions de front d'onde d'une repr\'esentation unitaire $\pi$ et de repr\'esentation fortement tra\c cable~: Le front d'onde de $\pi$ est une partie conique ferm\'ee du fibr\'e cotangent $T^*G$ priv\'e de la section nulle, invariante par translation \`a gauche et \`a droite, d\'efinie comme l'adh\'erence de la r\'eunion des fronts d'onde $WF(\mop{Tr}_\pi(T))$ o\`u $T$ parcourt l'ensemble des op\'erateurs \`a trace sur l'espace de $\pi$, et o\`u $\mop{Tr}_\pi(T)(g)=\mop{Tr}(\pi(g)T)$. Une repr\'esentation $\pi$ est fortement tra\c cable s'il existe un \'el\'ement $v$ formellement positif de l'alg\`ebre enveloppante $\cal U(\g g)$ tel que $\pi(v)$ soit essentiellement auto-adjoint et inversible sur l'espace des vecteurs $C^\infty$, et d'inverse \`a trace.
\ssq
Dans la premi\`ere partie nous donnons trois crit\`eres pour qu'une repr\'esentation unitaire $\pi$ soit fortement tra\c cable (proposition I.1). En particulier $\pi$ est fortement tra\c cable si et seulement s'il existe un entier $m$ tel que pour toute fonction $\varphi\in C^m_c(G)$ l'op\'erateur $\pi(\varphi)$ est \`a trace. L'existence de repr\'esentations tra\c cables non fortement tra\c cables est un probl\`eme ouvert (aimablement communiqu\'e par R. Howe).
\ssq
Dans la seconde partie nous rappelons quelques r\'esultats de R. Howe sur le front d'onde d'une repr\'esentation unitaire. Le lien entre ces deux notions est le suivant~: si $\pi$ est une repr\'esentation unitaire fortement tra\c cable, son front d'onde en l'\'el\'ement neutre co\"\i ncide avec le front d'onde de son caract\`ere-distribution [Hw theorem 1.8]. La d\'emonstration de ce r\'esultat \'etant donn\'ee par R. Howe dans le cas unimodulaire, nous montrons qu'elle reste valable dans le cas g\'en\'eral (th\'eor\`eme II.2).
\ssq
La troisi\`eme partie est consacr\'ee \`a la d\'emonstration du r\'esultat principal (th\'eor\`eme III.7)~: pour toute distribution $\varphi$ \`a support compact sur $G$ telle que~:
$$WF(\varphi)\cap -WF_\pi=\emptyset,$$
l'op\'erateur a priori non born\'e $\pi(\varphi)$ que l'on peut d\'efinir sur le domaine des vecteurs $C^\infty$ est en fait born\'e, et m\^eme r\'egularisant, c'est-\`a-dire qu'il envoie l'espace de la repr\'esentation dans l'espace des vecteurs $C^\infty$. Dans le cas o\`u $\pi$ est de plus fortement tra\c cable, l'op\'erateur $\pi(\varphi)$ est \`a trace.
\ssq
L'ingr\'edient essentiel de la d\'emonstration est le noyau de la chaleur $(p_t)_{t>0}$ sur le groupe, qui fournit une famille $\pi(p_t)$ croissante d'op\'erateurs born\'es qui converge fortement vers l'identit\'e. Pour tout op\'erateur de Hilbert-Schmidt $T$ sur l'espace de la repr\'esentation respectant le domaine de $\pi(\varphi)$ on montre que $\pi(\varphi)T$ est de Hilbert-Schmidt et que \hbox{$\pi(p_t*\varphi)T$} converge vers $\pi(\varphi)T$ en norme de Hilbert-Schmidt lorsque $t\rightarrow 0$. Plus pr\'ecis\'ement nous montrons que le front d'onde de $\varphi^**\varphi$ ne rencontre pas $WF_\pi$, ce qui nous permet de consid\'erer le produit des deux distributions $\mop{Tr}_\pi(TT^*)$ et $\varphi^**\varphi$, et nous d\'emontrons l'\'egalit\'e~:
$$\|\pi(\varphi)T\|_2^2=<\mop{Tr}_\pi(TT^*).(\varphi^**\varphi),\,1>.$$
\qquad Ceci entra\^\i ne le fait que $\pi(\varphi)$ est born\'e (corollaire III.5). Pour toute partie conique ferm\'ee $\Gamma$ de $T^*G-\{0\}$ ne rencontrant pas $WF_\pi$ on d\'esigne par $\Cal E'_\Gamma (G)$ l'espace des distributions \`a support compact dont le front d'onde est inclus dans $\Gamma$, muni d'une topologie localement convexe \`a base non d\'enombrable de voisinages ([Dui], [Ga]) dont nous rappelons la d\'efinition en III.1.
Nous ne pouvons pas montrer en g\'en\'eral la continuit\'e s\'equentielle de la correspondance $\varphi\rightarrow \pi(\varphi)$ de $\Cal E'_\Gamma (G)$ dans l'espace des op\'erateurs born\'es. Mais lorsque la repr\'esentation est fortement tra\c cable cette correspondance est s\'equentiellement continue de $\Cal E'_\Gamma (G)$ dans l'espace des op\'erateurs \`a trace.
\ssq
Dans la quatri\`eme partie nous montrons (th\'eor\`eme IV.1) que toute repr\'esentation unitaire irr\'eductible $\pi$ associ\'ee \`a une orbite coadjointe $\Omega$ par une ``bonne'' formule des caract\`eres de Kirillov est fortement tra\c cable. Plus pr\'ecis\'ement nous supposons que la repr\'esentation $\pi$ v\'erifie l'hypoth\`ese suivante ($j_G$ d\'esignant le jacobien de l'exponentielle~: cf notations 0.1)~:
\smallskip
{\bf Hypoth\`ese (H)\/}~: 
\par
{\sl
$\pi$ est associ\'ee \`a une orbite coadjointe {\sl ferm\'ee} $\Omega$ par la m\'ethode des orbites de Kirillov, de la mani\`ere suivante~: il existe un voisinage exponentiel $U$ de $0$ dans $\g g$, un entier strictement positif $\kappa(\pi)$ et une fonction $P_\Omega$ analytique partout non nulle sur $U$ et telle que $P(0)=1$, tels que la formule des caract\`eres~:
$$\mop{Tr}\pi(T)=\kappa(\pi)\int_\Omega \Cal F\bigl(P_\Omega\inver .j_G.(T\circ\exp)\bigr)(\omega)\, d\beta_\Omega(\omega)\leqno{(*)}$$
soit v\'erifi\'ee pour toute fonction $T$ de classe $C^\infty$ \`a support compact inclus dans $\exp U$.\/}
\smallskip
R. Howe a indiqu\'e dans le cas nilpotent comment calculer le front d'onde d'une telle repr\'esentation [Hw proposition 2.2]~:
$$WF_\pi(e)=-AC(\Omega)$$
o\`u $AC(\Omega)$ d\'esigne le c\^one asymptote \`a l'orbite. Nous donnons une d\'emonstration de ce r\'esultat et nous montrons (th\'eor\`eme IV.3.1)) que l'inclusion~: $WF_\pi(e)\subset -AC(\Omega)$ reste valable dans le cadre g\'en\'eral de l'hypoth\`ese (H). Cette derni\`ere inclusion est un cas particulier d'un r\'esultat g\'en\'eral sur le front d'onde de la transform\'ee de Fourier d'une mesure positive temp\'er\'ee [D-V lemme 4]. Nous ignorons si l'\'egalit\'e v\'erifi\'ee dans le cas nilpotent reste valable dans le cadre g\'en\'eral de l'hypoth\`ese (H).
\ssq
Utilisant le th\'eor\`eme III.7 nous montrons (th\'eor\`eme IV.3.2)) que la formule des caract\`eres est v\'erifi\'ee pour tous les op\'erateurs $\pi(\varphi)$ o\`u $\varphi$ est une distribution \`a support compact contenu dans un voisinage convenable de l'\'el\'ement neutre telle que~:
$$WF(\varphi)\cap -WF_\pi=\emptyset.$$
En particulier dans le cas des op\'erateurs pseudo-diff\'erentiels au sens de [M2], qui sont des op\'erateurs $\pi(\varphi)$ o\`u $\varphi$ est une distribution \`a support compact et \`a support singulier r\'eduit \`a $\{e\}$, la condition ci-dessus implique au vu du th\'eor\`eme IV.3~:
$$WF(\varphi)\cap AC(\Omega)=\emptyset.$$
Dans la cinqui\`eme partie nous comparons nos r\'esultats avec un r\'esultat ant\'erieur [M2 th\'eor\`eme III.1] dont la d\'emonstration \'etait incompl\`ete~: ce r\'esultat appara\^\i t ici comme une cons\'equence directe du th\'eor\`eme III.7, alors que la d\'emonstration succincte qui en est donn\'ee dans [M2] s'appuie quant \`a elle implicitement sur le th\'eor\`eme IV.3.2) du pr\'esent article. Nous donnons enfin une application de nos r\'esultats \`a une conjecture de M. Duflo et M. Vergne [D-V] sur la restriction d'une repr\'esentation $\pi$ fortement tra\c cable \`a un sous-groupe ferm\'e transverse au front d'onde de $\pi$.
\bigskip
{\bf Remerciements~:} Je remercie vivement Detlef M\"uller d'avoir attir\'e mon attention sur cette lacune dans la d\'emonstration du th\'eor\`eme III.1 de mon article ant\'erieur [M2], ainsi que pour ses remarques utiles concernant la troisi\`eme partie. Je remercie \'egalement Abderrazak Bouaziz pour d'utiles discussions concernant la derni\`ere application (paragraphe V.2).
\alinea{Notations et rappels~:}
(0.1). On d\'esignera par $G$ un groupe de Lie r\'eel connexe de dimension $n$, d'alg\`ebre de Lie $\g g$ et de dual $\g g^*$. On d\'esigne par $dx$ une mesure de Lebesgue sur $\g g$, et par $dg$ une mesure de Haar \`a gauche sur $G$ normalis\'ee de telle fa\c con que le jacobien $j_G$ de l'exponentielle s'\'ecrive~:
$$j_G(x)=\bigl|\det({1-e^{-\smop{ad}x} \over \mop{ad}x})\bigr|$$
On d\'esigne par $\Delta_G=\det\mop{Ad}g$ la fonction module. Par le choix de la mesure de Haar $dg$ nous identifierons l'espace $C^\infty(G)$ et l'espace des densit\'es $C^\infty$ sur $G$. Par dualit\'e l'espace des distributions sur $G$ s'identifie \`a l'espace des fonctions g\'en\'eralis\'ees sur $G$.
\smallskip
(0.2). Soit $U$ un voisinage de $0$ dans $\g g$ tel que l'exponentielle soit un diff\'eomorphisme de $U$ sur son image. Pour toute distribution $\varphi$ sur $G$ \`a support compact inclus dans $\exp U$ on d\'efinit la distribution $\wt \varphi$ par~:
$$\wt \varphi=j_G.\exp^* \varphi$$
(0.3). Soit $V$ un ouvert de $G$. Soit $m\in \N\cup\{+\infty\}.$ On d\'esigne par $C^m(V)$ (resp. $C^m_c(V)$) l'espace des fonctions (resp. des fonctions \`a support compact) $m$ fois diff\'erentiables sur $V$.
\smallskip
(0.4). On d\'esigne par $\Cal D'(V)$ l'espace des distributions sur $V$, et par $\Cal E'(V)$ l'espace des distributions \`a support compact sur $V$. On met sur $\Cal E'(V)$ la topologie limite inductive des $\Cal D'(K)$ o\`u $K$ parcourt l'ensemble des compacts de $V$. En particulier une suite $\varphi_k$ converge dans $\Cal E'(V)$ si et seulement si elle converge dans $\Cal D'(V)$ et tous les supports des $\varphi_k$ sont contenus dans un m\^eme compact $K$.
\smallskip
(0.5). On d\'esignera par $T^*G\backslash \{0\}$ le fibr\'e cotangent de $G$ priv\'e de la section nulle. Une partie $C$ de $T^*G\backslash\{0\}$ est {\sl conique\/} si pour tout $(g,\xi)\in C$, $(g,t\xi) \in C$ pour tout r\'eel $t>0$. On notera alors $-C$ l'ensemble des $\{(g,-\xi),\, (g,\xi)\in C\}$.
\smallskip
(0.6). Soit $\Cal H$ un espace de Hilbert, qui sera toujours suppos\'e s\'eparable. On d\'esigne par $J_p$ le $p^{\hbox{\sevenrm i\`eme}}$ id\'eal de Schatten de $\Cal H$. En particulier $J_\infty$ d\'esigne l'espace des op\'erateurs born\'es sur $\Cal H$, $J_1$ et $J_2$ les id\'eaux de $J_\infty$ form\'es des op\'erateurs \`a trace et des op\'erateurs de Hilbert-Schmidt respectivement.
\smallskip
(0.7). On d\'esignera par $\pi$ une repr\'esentation unitaire fortement continue du groupe de Lie $G$ dans un espace de Hilbert s\'eparable qui sera not\'e $\Cal H_\pi$. On d\'esignera alors par $\Cal H^\infty_\pi$ l'espace des vecteurs ind\'efiniment diff\'erentiables de la repr\'esentation. L'espace $\Cal H_\pi^\infty$ est constitu\'e des vecteurs $u$ tels que pour tout $v\in\Cal H_\pi$ le coefficient~:
$$C_{u,v}:g\longmapsto <\pi(g)u,\,v>$$
soit $C^\infty$ sur $G$.
\smallskip
(0.8). On d\'efinit pour tout $\varphi\in\Cal E'(G)$ l'op\'erateur (en g\'en\'eral non born\'e) $\pi(\varphi)$ de domaine $\Cal H_\pi^\infty$ par la formule~:
$$<\pi(\varphi)u,\, v>=<\varphi,\, C_{u,v}>$$ 
(Voir [J\o] pour une d\'efinition dans le cadre g\'en\'eral des repr\'esentations dans un espace de Banach). La correspondance $\varphi\mapsto \pi(\varphi)$ est continue de $\Cal E'(G)$ vers l'espace des op\'erateurs de domaine $\Cal H_\pi^\infty$ muni de la topologie de la convergence faible g\'en\'eralis\'ee au sens de Kato.
\smallskip
(0.9). Pour toute distribution $\varphi\in \Cal E'(G)$ on pose~:
$$\varphi^*=\Delta_G.\overline {i^*\varphi}$$
o\`u $\Delta_G$ d\'esigne la fonction module et $i$ le diff\'eomorphisme $g\mapsto g\inver$. On a pour tout couple $(u,v)$ dans $\Cal H_\pi^\infty$~:
$$<\pi(\varphi)u,\,v>=<u,\pi(\varphi^*)v>.$$
\prop{0.1}
Pour tout couple $(S,T)$ de distributions \`a support compact sur $G$ et pour tout $u$ vecteur $C^\infty$ de la repr\'esentation on a~:
$$\pi(S*T)u=\pi(S)\pi(T)u$$
\dem
C'est clair dans le cas o\`u $S,T$ sont dans $C^\infty_0(G)$. Soit $U$ un voisinage exponentiel de $0$ dans $\g g$. On se donne une fonction $\alpha\in C^\infty_0(\exp U)$ positive d'int\'egrale $1$, et on consid\`ere l'approximation de l'identit\'e~:
$$\alpha_k(\exp x)=k^n\alpha(\exp (kx))$$
Posons $T_k=\alpha_k*T$ et $S_k=\alpha_k*S$.La distribution $S_k*T_k$ converge vers $S*T$ pour la topologie de $\Cal E'(G)$. Pour tout $u\in\Cal H_\pi^\infty$ et pour tout $v$ le scalaire $<\pi(S_k)\pi(T_k)u,\, v>$ converge donc vers $<\pi(S)\pi(T)u,\, v>$, et $<\pi(S_k*T_k)u,\,v>$ converge vers $<\pi(S*T)u,\,v>$. La proposition d\'ecoule alors de l'\'egalit\'e~:
$$\pi(S_k*T_k)=\pi(S_k)\circ \pi(T_k)$$
\qed
\paragraphe{I. Repr\'esentations fortement tra\c cables}
\qquad Soit $G$ un groupe de Lie, $\g g$ son alg\`ebre de Lie et $\Cal U(\g g)$ son alg\`ebre envelopppante, que nous identifierons avec l'alg\`ebre des distributions de support $\{e\}$ sur $G$. On d\'esigne par $v\mapsto v^*$ l'unique anti-automorphisme de $\Cal U(\g g)$ tel que $x^*=-x$ pour $x\in\g g$. Un \'el\'ement $v\in\Cal U(\g g)$ est dit {\sl formellement positif\/} [Hw \S\ I] s'il s'\'ecrit sous la forme~:
$$u=\sum_{i=1}^k u_i^*u_i$$
Un exemple d'\'el\'ement formellement positif est le laplacien~:
$$\Delta=-\sum_1^n X_i^2$$
pour une base $(\uple Xn)$ de $\g g$.
Une repr\'esentation unitaire $\pi$ de $G$ est dite {\sl fortement tra\c cable\/} [Hw \S\ I] s'il existe un \'el\'ement $u\in\Cal U(\g g)$ formellement positif tel que $\pi(u)$ soit inversible sur $\Cal H_\pi^\infty$ et d'inverse \`a trace.
\prop{I.1}
Soit $\pi$ une repr\'esentation unitaire d'un groupe de Lie $G$. Les quatre conditions suivantes sont \'equivalentes~:
\smallskip
1) $\pi$ est fortement tra\c cable
\smallskip
2) Il existe un entier $s_0$ tel que pour tout entier $s\ge s_0$ l'op\'erateur $\pi(1+\Delta)^{-s}$ est \`a trace.
\smallskip
3) Il existe un voisinage $V$ de $e$ dans $G$ et un entier $m$ tel que pour toute fonction $\varphi\in C^m_c(V)$ l'op\'erateur $\pi(\varphi)$ est \`a trace.
\smallskip
4) Il existe un entier $m$ tel que pour tout $\varphi\in C^m_c(G)$ l'op\'erateur $\pi(\varphi)$ est \`a trace.
\dem
On montre facilement que si $u$ est formellement positif il en est de m\^eme de $u^k$ pour tout $k$. 2)$\Longrightarrow$ 1) est alors imm\'ediat. L'implication 1)$\Longrightarrow$ 4) est due \`a R. Howe~: supposons que $\pi$ soit fortement tra\c cable, et soit $v\in \Cal U(\g g)$ formellement positif tel que $\pi(v)$ soit inversible d'inverse \`a trace. Alors si $\varphi\in C^m_c(G)$ avec $m$ assez grand la convol\'ee $v*\varphi$ est une fonction continue \`a support compact sur $G$. On d\'eduit 4) de l'\'egalit\'e entre op\'erateurs born\'es~:
$$\pi(\varphi)=\pi(v)\inver\pi(v*\varphi).$$
4)$\Longrightarrow$ 3) \'etant imm\'ediat, reste l'implication 3)$\Longrightarrow$ 2). On suppose que 3) est v\'erifi\'ee, et on restreint s'il le faut le voisinage $V$ de fa\c con \`a ce que l'op\'erateur de convolution par $(1+\Delta)^s$ soit elliptique sur $V$. Il existe alors deux fonctions $\varphi$ et $\psi$ dans $C^{2s-n-1}_c(V)$ et $C^\infty_c(V)$ respectivement, telles que~:
$$(1+\Delta)^s *\varphi=\delta_0+\psi.$$
L'op\'erateur $\pi(1+\Delta)$ est essentiellement auto-adjoint [Ne-St th. 2.2], et son spectre est contenu dans $[1,+\infty[$. D'apr\`es le th\'eor\`eme spectral cet op\'erateur est inversible sur son domaine et son inverse est born\'e. On peut donc \'ecrire~:
$$\pi(1+\Delta)^{-s}=\pi(\varphi)-\pi(1+\Delta)^{-s}\pi(\psi)$$
Si l'entier $s$ est assez grand pour que l'on ait $2s-n-1\ge m$ il est donc clair que $\pi(1+\Delta)^{-s}$ est \`a trace.
\qed
{\sl Remarque\/}~: consid\'erons la famille d'espaces de Sobolev $\Cal H_\pi^t$ introduits par R. Goodman [Go] d\'efinis par compl\'etion de $\Cal H_\pi^\infty$ pour les normes~:
$$\|u\|_t=\|\pi(1+\Delta)^{\frac t2}u\|$$ 
On voit que la condition 2) de la proposition entra\^\i ne que la famille des espaces $\Cal H_\pi^t$ est nucl\'eaire, c'est-\`a-dire que l'inclusion de $\Cal H_\pi^t$ dans $\Cal H_\pi^{t-s}$ est \`a trace pour $s$ assez grand. Le pas $s$ peut ici \^etre choisi ind\'ependant de $t$.
\goodbreak
\paragraphe{II. Front d'onde d'une repr\'esentation}
\qquad Dans ce paragraphe nous rappelons la notion de front d'onde d'une repr\'esentation unitaire et nous donnons la d\'emonstration du fait que le front d'onde d'une repr\'esentation fortement tra\c cable et le front d'onde de son caract\`ere co\"\i ncident en l'\'el\'ement neutre. Nous reprenons simplement la d\'emonstration que R. Howe donne dans le cas unimodulaire et nous montrons qu'elle est valable en g\'en\'eral. Nous avons choisi d'enlever la section nulle du front d'onde, comme dans le cas du front d'onde d'une distribution. Soit $\pi$ une repr\'esentation unitaire d'un groupe de Lie $G$. Pour tout $T\in J_1(\Cal H_\pi)$ on consid\`ere la fonction~:
$$\mop{Tr}_\pi(T):g\longrightarrow \mop{Tr}\bigl(\pi(g)T\bigr).$$
que l'on peut aussi voir comme une distribution sur $G$. Le front d'onde $WF_\pi$ est alors d\'efini comme l'adh\'erence dans $T^*G\backslash \{0\}$ de la r\'eunion~:
$$\bigcup_{T\in J_1}WF\bigl(\mop{Tr}_\pi(T)\bigr).$$
\qquad
L'ensemble $WF_\pi$ est invariant par translation \`a gauche et \`a droite. Il est donc enti\`erement d\'etermin\'e par son intersection $WF^0_\pi$ avec $T^*_eG$, qui est un c\^one ferm\'e $\mop{Ad}^*G$-invariant dans $\g g^*$ priv\'e de $\{0\}$. Le front d'onde d'une repr\'esentation est caract\'eris\'e par les propri\'et\'es suivantes~:
\prop{II.1 \rm(R. Howe)}
Soit $\pi$ une repr\'esentation unitaire d'un groupe de Lie $G$ d'alg\`ebre de Lie $\g g$, soit $U$ un ouvert de $\g g^*$ et soit $\Cal U$ l'ouvert invariant par translation \`a gauche de $T^*G$ tel que $\Cal U\cap T^*_eG=U$. Alors les conditions suivantes sont \'equivalentes~:
\smallskip
1) $U\cap WF_\pi^0=\emptyset$.
\smallskip
2) Pour tout $T\in J_1$, pour tout voisinage $V$ de $e$ dans $G$ et pour tout couple $(\varphi_\alpha, \psi_\alpha)$ d\'ependant de mani\`ere $C^\infty$ d'un param\`etre $\alpha\in\R^k$ o\`u $\varphi_\alpha\in C^\infty_c(V)$, et $\psi_\alpha\in C^\infty(G)$ \`a valeurs r\'eelles telles que $d\psi_\alpha(\mop{supp}\varphi_\alpha)\subset\Cal U$, l'int\'egrale~:
$$I(\varphi_\alpha,\psi_\alpha,T)(t)=\int_G \mop{Tr}_\pi(T)(g)\varphi_\alpha(g)e^{it\psi_\alpha(g)}\,dg$$
est \`a d\'ecroissance rapide lorque $t\rightarrow +\infty$, et on a les estimations~:
$$\mopl{sup}_{t\ge 1}t^N|I(\varphi_\alpha,\psi_\alpha,T)(t)|\le C_N(\varphi_\alpha,\psi_\alpha)\|T\|_1$$
les constantes $C_N$ pouvant \^etre choisies ind\'ependamment de $\alpha$ lorsque $\alpha$ varie dans un compact de $\R^k$.
\smallskip
3) pour tout voisinage $V$ de $e$ dans $G$, pour tout couple $(\varphi_\alpha, \psi_\alpha)$ d\'ependant de mani\`ere $C^\infty$ d'un param\`etre $\alpha\in\R^k$ o\`u $\varphi_\alpha\in C^\infty_c(V)$, et $\psi_\alpha\in C^\infty(G)$ \`a valeurs r\'eelles telles que $d\psi_\alpha(\mop{supp}\varphi_\alpha)\subset\Cal U$, la norme d'op\'erateur de $\pi(\varphi_\alpha e^{it\psi_\alpha})$ est \`a d\'ecroissance rapide lorque $t\rightarrow +\infty$, et les semi-normes~:
$$\gamma_N(\varphi_\alpha,\psi_\alpha)=\mopl{sup}_{t\ge 1}t^N\|\pi(\varphi_\alpha e^{it\psi_\alpha})\|_\infty$$
sont born\'ees uniform\'ement lorque $\alpha$ varie dans un compact de $\R^k$.
\dem
d'apr\`es l'\'egalit\'e~:
$$I(\varphi_\alpha,\psi_\alpha,T)(t)=\mop{Tr}\bigl(\pi(\varphi_\alpha e^{it\psi_\alpha})T\bigr)$$
et la dualit\'e entre $J_1$ et $J_\infty$~:
$$\forall A\in J_\infty, \|A\|_\infty=\mopl{sup}_{T\in J_1-\{0\}}
{|\mop{Tr}AT|\over \|T\|_1},$$
les conditions 2) et 3) sont \'equivalentes, et dans la condition 2) les meilleures constantes sont~:
$$ C_N(\varphi_\alpha,\psi_\alpha)= \gamma_N(\varphi_\alpha,\psi_\alpha).$$
\qquad
Par ailleurs 1) implique (par bi-invariance du front d'onde $WF_\pi$) que $\Cal U$ ne rencontre pas $WF_\pi$. Ceci entra\^\i ne 2) de par la d\'efinition du front d'onde de $\mop{Tr}_\pi(T)$ [Dui proposition 1.3.2]. R\'eciproquement 2) implique que $\Cal U$ est disjoint du front d'onde de $\mop{Tr}_\pi(T)$ pour tout $T\in J_1$, ce qui entra\^\i ne 1) par bi-invariance de $WF_\pi$.
\qed
\th{II.2 \rm (R.Howe)}
Soit $\pi$ une repr\'esentation unitaire fortement tra\c cable d'un groupe de Lie $G$. Soit $\chi_\pi$ le caract\`ere distribution de $\pi$, et soit $WF^0(\chi_\pi)$ l'intersection du front d'onde de $\chi_\pi$ avec $T^*_e(G)$. Alors on a~:
$$\eqalign{WF(\chi_\pi)	&\subset WF_\pi	\cr
	WF_\pi^0	&=WF^0(\chi_\pi).\cr}$$
\dem
Pour montrer la premi\`ere inclusion, on consid\`ere un \'el\'ement $v\in\Cal U(\g g)$ formellement positif tel que $\pi(v)$ soit inversible et d'inverse $T$ dans $J_1$. Pour tout $\varphi\in C^\infty_c(G)$ on a~:
$$\eqalign{\chi_\pi(\varphi)	&=\mop{Tr}\bigl(\pi(\varphi)\pi(v)T \bigr)\cr
				&=\mop{Tr}\bigl(\pi(\varphi *v)T \bigr)\cr
				&=\mop{Tr}_\pi(T)(\varphi *v),\cr}$$
c'est-\`a-dire~:
$$\chi_\pi=\mop{Tr}_\pi(T)*v^*.$$
Comme la convolution par $v^*$ respecte le front d'onde on a les inclusions~:
$$WF(\chi_\pi)\subset WF(\mop{Tr}_\pi(T))\subset WF_\pi.$$
\qquad
R\'eciproquement soit $\xi$ non nul dans $\g g^*\backslash WF^0(\chi_\pi)$, soit $U$ un voisinage relativement compact de $\xi$ disjoint de $WF^0(\chi_\pi)$, et soit $\Cal U$ le sous-ensemble de $T^*G$ construit \`a partir de $U$ par translation \`a gauche. Comme $WF(\chi_\pi)$ est ferm\'e il existe un voisinage $V$ de $e$ dans $G$ tel que $\Cal U\cap T^*V$ soit disjoint de $WF(\chi_\pi)$.
\ssq
Soit $V_1$ un voisinage sym\'etrique de $e$ dans $G$ tel que $V_1^2\subset V$, que l'on supposera relativement compact, et soit $\varphi\in C^\infty_c(V_1)$, $\psi\in C^\infty(V)$ telles que $d\psi(\mop{supp}\varphi)\subset \Cal U$. Alors $\chi_\pi\bigl(L_g.(\varphi e^{it\psi}) \bigr)$ est \`a d\'ecroissance rapide lorsque $t\rightarrow +\infty$, et ceci uniform\'ement pour tout $g\in V_1$.
\ssq
On pose $\varphi_t=\varphi e^{it\psi}$, et on calcule explicitement la norme de Hilbert-Schmidt de $\pi(\varphi_t)$, ce qui nous donne~:
$$\|\pi(\varphi_t)\|_2^2=\int_G \varphi_t^*(y)\chi_\pi\bigl(L_y\varphi_t \bigr)\, dy$$
qui est donc \`a d\'ecroissance rapide d'apr\`es ce qui pr\'ec\`ede (ici $dy$ d\'esigne une mesure de Haar \`a gauche sur $G$). C'est vrai aussi pour la norme $\|\pi(\varphi_t)\|_\infty$, et il est clair que si les fonctions $\varphi$ et $\psi$ d\'ependent de mani\`ere $C^\infty$ d'un param\`etre auxiliaire $\alpha$ variant dans un compact de $\R^k$ les semi-normes~:
$$\mopl{sup}_{t\ge 1}t^N\|\pi(\varphi_t)\|_\infty$$
sont born\'ees ind\'ependamment de $\alpha$. Le crit\`ere 3) de la proposition II.1 dit alors que $U$ est disjoint de $WF^0_\pi$, ce qui d\'emontre le th\'eor\`eme.
\qed
{\sl Remarque\/}~: l'hypoth\`ese de tra\c cabilit\'e forte n'intervient que dans la d\'emonstration de la premi\`ere inclusion. Il est donc vrai que pour toute repr\'esentation {\sl tra\c cable\/} on a~:
$$WF_\pi^0\subset WF^0(\chi_\pi).$$
\goodbreak
\paragraphe{III. Une classe d'op\'erateurs r\'egularisants}
\alinea{III.1. Quelques espaces de distributions}

\qquad Nous rappelons ici la d\'efinition et quelques propri\'et\'es des espaces $\Cal D'_\Gamma(X)$ introduits par L. H\"ormander ([Hr \S 2.5], [Ga], [Dui]), ainsi que de leurs variantes \`a support compact $\Cal E'_\Gamma(X)$.  
\ssq 
Soit $X$ une vari\'et\'e, soit $\Gamma$ une partie conique ferm\'ee de l'espace cotangent $T^*X$ priv\'e de la section nulle. On notera pour tout $x\in X$~:
$$\Gamma_x=\Gamma\cap T^*_xX.$$
On consid\`ere l'espace $\Cal D'_\Gamma(X)$ des distributions \`a support compact sur $X$ dont le front d'onde est contenu dans $\Gamma$. On met sur $\Cal D'_\Gamma(X)$ la topologie (localement convexe \`a base non d\'enombrable de voisinages) d\'efinie par les semi-normes de la topologie de $\Cal D'(X)$ et les semi-normes~:
$$C_{N,\alpha,\psi}(\varphi)=\mopl{sup}_{\tau\ge 1}\tau^N|<e^{-i\tau \psi}\alpha,\varphi>|$$
o\`u $N$ est un entier positif, $\alpha\in C^\infty_c(X)$ et $\psi\in C^\infty(X)$ est une fonction \`a valeurs r\'eelles telle que $d\psi(g)\notin \Gamma$ pour tout $g\in \mop{supp}u$ [Dui\S\ 1.3]. On rappelle [Dui theorem 1.3.6] que si $\Gamma_1$ et $\Gamma_2$ sont deux parties coniques ferm\'ees de l'espace cotangent $T^*X$ priv\'e de la section nulle telles que $\Gamma_1\cap -\Gamma_2=0$ le produit des fonctions $C^\infty$ s'\'etend en une unique op\'eration bilin\'eaire continue~:
$$\eqalign{\Cal D'_{\Gamma_1}(X)\times \Cal D'_{\Gamma_1}(X)
		&\longrightarrow \Cal D'_{\overline{\Gamma_1+\Gamma_2}}(X)	\cr 
		(u,v)	&\longmapsto uv.\cr}$$
On notera~: $\Cal E'_\Gamma(X)$ l'ensemble des distributions \`a support compact appartenant \`a $\Cal D'_\Gamma(X)$. On met sur $\Cal E'_\Gamma(X)$ la topologie limite inductive des $\Cal D'_\Gamma(K)$ o\`u $K$ parcourt les compacts de $X$. C'est un espace localement convexe \`a base non d\'enombrable de voisinages, et il est clair que le produit d\'efini ci-dessus est continu de $\Cal E'_{ \Gamma_1}(X)\times \Cal E'_{\Gamma_2}(X)$ dans $\Cal E'_{\overline{\Gamma_1+\Gamma_2 }}(X)$.
\prop{III.1 {\rm [Dui proposition 1.3.4]}}
Soit $f$ une application $C^\infty$ d'une vari\'et\'e $X$ vers une vari\'et\'e $Y$. Pour toute partie conique ferm\'ee $\Gamma$ de $T^*X\backslash\{0\}$ on consid\`ere~:
$$f_*\Gamma=\{(y,\eta)\in T^*Y\backslash\{0\}/\exists x\in X,y=f(x)\hbox{ et } (x,f^*\eta)\in \Gamma\}.$$
Alors l'application $f$ se prolonge en une application lin\'eaire continue~:
$$f_*:\Cal E'_\Gamma(X)\longrightarrow\Cal E'_{\overline{f_*\Gamma}}(Y),$$
et si l'application $f$ est propre elle se prolonge en une application lin\'eaire continue~:
$$f_*:\Cal D'_\Gamma(X)\longrightarrow\Cal D'_{\overline{f_*\Gamma}}(Y).$$
\dem
La topologie de $\Cal D'_\Gamma(X)$ est d\'efinie par les semi-normes~:
$$\eqalign{C_u(\varphi)		&=|<\varphi,\,u>|\cr
C_{N,\alpha,\psi}(\varphi)	&=\mopl{sup}_{t\ge 1}t^N
				|<\varphi,\, \alpha e^{-it\psi}>|\cr}$$
pour $u,\alpha\in C^\infty_c(X)$, $N\in\N$ et $\psi\in C^\infty(X)$ telle que $d\psi(\mop{supp}\alpha)\cap\Gamma=\emptyset$. Si $f$ est propre $u\circ f$ et $\alpha\circ f$ sont dans $C^\infty_c(X)$, et on a \'egalit\'e entre semi-normes~:
$$C_u(f_*\varphi)=C_{u\circ f}(\varphi)\hbox{,\qquad}C_{N,\alpha,\psi}(f_*u)
=C_{N,\alpha\circ f,\psi\circ f}(u)$$
o\`u cette fois-ci $u,\alpha\in C^\infty_c(Y)$, $N\in\N$ et $\psi\in C^\infty(Y)$ telle que $d\psi(\mop{supp}\alpha)\cap f_*\Gamma=\emptyset$. La composition des diff\'erentielles~:
$$d_x(\psi\circ f)=f^*(d_{f(x)}\psi)$$
et la d\'efinition de $f_*\Gamma$ permettent de conclure \`a la continuit\'e de
$f_*:\Cal D'_\Gamma(X)\longrightarrow\Cal D'_{\overline{f_*\Gamma}}(Y)$. Le passage par limite inductive au cas des distributions \`a support compact est alors imm\'ediat.
\alinea{III.2. Front d'onde et convolution}
\qquad Nous montrons ici un r\'esultat sur le front d'onde de la convol\'ee de deux distributions sur un groupe de Lie $G$ dont l'une est \`a support compact. Ce r\'esultat d\'ecoule directement de la proposition III.1. 
\prop{III.2}
Soient $\varphi$ et $\psi$ deux distributions sur un groupe de Lie $G$ dont l'une est \`a support compact, et soit $z\in G$. Alors le front d'onde de $\varphi*\psi$ en $z$ v\'erifie~:
$$WF_z(\varphi*\psi)\subset\overline{\bigcup_{xy=z}
R_{y\inver}^*WF_x(\varphi)\cap L_{x\inver}^*WF_y(\psi)}$$
o\`u l'adh\'erence est prise dans $T^*_zG-\{0\}$. De plus si $\Gamma'$ et $\Gamma''$ sont deux parties coniques ferm\'ees de $T^*G\backslash\{0\}$ et si on consid\`ere la partie conique ferm\'ee $\Gamma$ de $T^*G\backslash\{0\}$ d\'efinie par~:
$$\Gamma_z=\overline{\bigcup_{xy=z}R_{y\inver}^*\Gamma'_x\cap L^*_{x\inver}\Gamma''_y}$$
la correspondance $(\varphi,\psi)\mapsto \varphi*\psi$~:
$$\Cal E'_{\Gamma'}(G)\times \Cal E'_{\Gamma''}(G)\longrightarrow
\Cal E'_\Gamma (G)$$
est s\'equentiellement continue.
\dem
l'application produit~:
$$m:G\times G\longrightarrow G$$
se prolonge \`a l'espace $\Cal D'_m(G\times G)$ des distributions $u$ sur $G\times G$ telles que la restriction de $m$ au support de $u$ est propre~:
$$m_*:\Cal D'_m(G\times G)\longrightarrow \Cal D'(G).$$
En particulier si $\varphi$ ou $\psi$ est \`a support compact sur $G$, le produit tensoriel $\varphi\otimes \psi$ appartient \`a $\Cal D'_m(G\times G)$ et on a par d\'efinition du produit de convolution~:
$$\varphi*\psi=m_*(\varphi\otimes \psi).$$
Soient $\Gamma$, $\Gamma'$ et $\Gamma''$ comme dans l'\'enonc\'e de la proposition. En appliquant la proposition III.1 on voit que la correspondance~: 
$$\eqalign{\Cal E'_{\Gamma'\times \Gamma''}(G\times G)	&\longrightarrow
\Cal E'_{\overline{m_*(\Gamma'\times \Gamma'')}}(G)\cr
\varphi\otimes\psi					&\longmapsto
\varphi*\psi\cr}$$ 
est continue. Par ailleurs la correspondance bilin\'eaire~:
$$\eqalign{\Cal D'_{\Gamma'(G)}\times \Cal D'_{\Gamma''(G)}
&	\longrightarrow \Cal D'_{\overline{\Gamma'\times \Gamma''}}(G\times G)\cr
(\varphi,\psi)	&\longmapsto \varphi\otimes \psi \cr}$$
est s\'equentiellement continue [Ga], et on voit par passage \`a la limite inductive sur les compacts que c'est \'egalement vrai en rempla\c cant $\Cal D'$ par $\Cal E'$. Il suffit donc de montrer que~:
$$\overline{m_*(\Gamma'\times \Gamma'')}\subset\Gamma.$$
La diff\'erentielle du produit s'\'ecrit~:
$$Dm(x,y;X,Y)=(xy,(R_y)_*X+(L_x)_*Y).$$
Sa transpos\'ee, de $T^*_{xy}G$ vers $T^*_{(x,y)}(G\times G)$ s'\'ecrit donc~:
$$m^*(xy,\eta)=(x,y;(R_y^*\otimes L_x^*)\delta\eta),$$
o\`u $\delta$ est la transpos\'ee de l'op\'eration d'addition dans $T_{xy}G$. On voit donc que~:
$$\eqalign{(m^*)\inver (\Gamma'_x\times\Gamma''_y)	&=\{(xy,\eta)/\delta\eta
\in R_{y\inver}^*\Gamma'_x\times R_{x\inver}^*\Gamma''_y \}\cr
	&\subset\{(xy,\eta)/\eta\in R_{y\inver}^*\Gamma'_x\cap R_{x\inver}^*\Gamma''_y\}.\cr}$$
Donc $\overline{m_*(\Gamma'\times \Gamma'')}\subset\Gamma$, ce qui montre la proposition.
\qed
Par bi-invariance du front d'onde d'une repr\'esentation unitaire on obtient imm\'ediatement le r\'esultat suivant~:
\cor{III.3}
1) Soit $\pi$ une repr\'esentation unitaire du groupe de Lie $G$, et soient $\varphi$ et $\psi$ deux distributions dans $\Cal E'(G)$ telles que $WF(\varphi)$ et $WF(\psi)$ ne rencontrent pas $-WF_\pi$. Alors $WF(\varphi*\psi)$ ne rencontre pas non plus $-WF_\pi$.
\smallskip
2) Soient $\Gamma_1$ et $\Gamma_2$ deux parties coniques ferm\'ees de $T^*G\backslash\{0\}$, et soit $\Cal G$ une partie conique ferm\'ee bi-invariante de $T^*G\backslash\{0\}$ contenant $\Gamma_1$ et $\Gamma_2$. Alors le produit de convolution est s\'equentiellement continu de $\Cal E'_{\Gamma_1}(G)\times \Cal E'_{\Gamma_2}(G)$ dans $\Cal E'_{\Cal G}(G)$.
\ndem
\alinea{III.3. Application \`a l'\'etude des op\'erateurs $\pi(\varphi)$}
\qquad Nous avons d\'efini dans l'introduction l'op\'erateur non born\'e $\pi(\varphi)$ pour une repr\'esentation unitaire $\pi$ d'un groupe de Lie $G$ et une distribution \`a support compact $\varphi$ sur $G$. Nous donnons ici une condition suffisante, portant sur le front d'onde de $\varphi$, pour que $\pi(\varphi)$ soit r\'egularisant, et nous montrons que dans le cas o\`u la repr\'esentation $\pi$ est fortement tra\c cable, la m\^eme condition implique que $\pi(\varphi)$ est un op\'erateur \`a trace.
\ssq
\th{III.4}
Soit $\pi$ une repr\'esentation unitaire d'un groupe de Lie $G$ et $\Gamma$ une partie conique ferm\'ee de l'espace cotangent $T^*G$ priv\'e de la section nulle ne rencontrant pas $-WF_\pi$. Alors tout op\'erateur $T\in J_2(\Cal H_\pi)$ laissant le domaine de $\pi(\varphi)$ invariant d\'efinit une correspondance lin\'eaire s\'equentiellement continue~:
$$\eqalign{R_T:\Cal E'_\Gamma(G)	&\longrightarrow J_2	\cr
		\varphi			&\longmapsto \pi(\varphi)T\cr}$$
\dem
On se donne une distribution $\varphi\in\Cal E'_\Gamma(G)$ et on consid\`ere le laplacien~:
$$\Delta=-\sum_{i=1}^n X_i$$ 
o\`u $(\uple Xn)$ est une base de $\g g$, et on consid\`ere le noyau de la chaleur $p_t=e^{-t\Delta}\delta$ associ\'e ([Ne], [Hu], [Rob],...). C'est une famille $(p_t)_{t>0}$ de fonctions $C^\infty$ sur $G$ qui v\'erifie~:
$$\eqalign{	&1)\ p_t^*=p_t	\cr
		&2)\ \|p_t\|_1=1	\cr
		&3)\ p_t(g)>0 \hbox{  pour tout $g\in G$}	\cr
		&4)p_t\fleche 8_{t\rightarrow 0} \delta
					\hbox{  dans $\Cal D'(G)$}\cr
		&5)p_t*p_s=p_{t+s}\hbox {  (propri\'et\'e de semi-groupe)}.\cr}$$
Soit $\pi$ une repr\'esentation unitaire de $G$. Alors [Ne \S \ 8 theorem 4] l'op\'erateur $\pi(p_t)$ est r\'egularisant. On sait aussi [Ne-St theorem 2.2] que $\pi(\Delta)$ est essentiellement auto-adjoint (et positif). On peut donc \'ecrire, gr\^ace au th\'eor\`eme spectral~:
$$\eqalign{\pi(\Delta)	&=\int_0^{+\infty}\lambda\,dE(\lambda)	\cr
	\pi(p_t)	&=\int_0^{+\infty}e^{-t\lambda}\,dE(\lambda)	\cr}$$
On voit ainsi que $(\pi(p_t))_{t\rightarrow 0}$ est une famille croissante d'op\'erateurs born\'es qui converge fortement vers l'identit\'e. On calcule alors la norme de Hilbert-Schmidt de $\pi(p_t)\pi(\varphi)T$ pour tout $T\in J_2$ laissant invariant le domaine de $\pi(\varphi)$~:
$$\eqalign{\|\pi(p_t)\pi(\varphi)T \|_2^2	&=\mop{Tr}\bigl(
\pi(p_t)\pi(\varphi)TT^*\pi(\varphi^*)\pi(p_t) \bigr)\cr
						&=\mop{Tr}_\pi(TT^*)(\varphi^**p_{2t}*\varphi).\cr}$$
Par bi-invariance du front d'onde de la repr\'esentation, $-WF_\pi\cap WF(\varphi^*)=\emptyset$. D'apr\`es le corollaire III.3 le front d'onde de $\varphi^**\varphi$ ne rencontre pas $-WF_\pi$. D'apr\`es l'hypoth\`ese du th\'eor\`eme il est donc licite d'effectuer le produit de la distribution $\varphi^**\varphi$ par la distribution $\mop{Tr}_\pi(TT^*)$, et on a~:
$$\mop{Tr}_\pi(TT^*)(\varphi^**p_{2t}*\varphi)\fleche 8_{t\rightarrow 0}
<\mop{Tr}_\pi(TT^*).(\varphi^**\varphi),\, 1>.$$
Autrement dit $\|\pi(p_t)\pi(\varphi)T \|_2^2$ admet une limite $\ell$ lorsque $t$ tend vers $0$. D'autre part, si $(e_i)$ d\'esigne une base orthonorm\'ee de $\Cal H_\pi$ contenue dans le domaine de $\pi(\varphi)$ on a~:
$$\|\pi(p_t)\pi(\varphi)T \|_2^2=\sum_{i\ge 0}<\pi(p_{2t})\pi(\varphi)Te_i,\,
\pi(\varphi)Te_i)>.$$
On en d\'eduit d'apr\`es le th\'eor\`eme de convergence monotone (puisque $\pi(p_{2t})$ est croissante) que la somme~:
$$\sum_{i\ge 0}\|\pi(\varphi)Te_i\|^2$$
converge vers la limite $\ell$ d\'efinie ci-dessus. Il r\'esulte alors de [DS lemma XI.9.32] que l'op\'erateur $\pi(\varphi)T$ est dans $J_2$ et que~:
$$\|\pi(\varphi)T\|_2^2=\ell=<\mop{Tr}_\pi(TT^*).(\varphi^**\varphi),\, 1>.$$
Enfin soit $\Cal G$ une partie conique ferm\'ee bi-invariante de $T^*G\backslash\{0\}$ contenant $\Gamma$ et ne rencontrant pas $-WF_\pi$. D'apr\`es le corollaire III.3 la correspondance $\varphi\mapsto \varphi^**\varphi$ est s\'equentiellement continue de $\Cal E'_\Gamma(G)$ dans $\Cal E'_{\Cal G}(G)$. La continuit\'e s\'equentielle de $R_T:\varphi\mapsto\pi(\varphi)T$ de $\Cal E'_\Gamma(G)$ dans $J_2$ provient alors directement de la continuit\'e du produit de $ \Cal E'_{WF_\pi}(G)\times\Cal E'_{\Cal G}(G)$ dans $\Cal E'_{\overline{\Cal G+WF_\pi}}(G)$.
\qed 
\cor{III.5}
Soit $\pi$ une repr\'esentation unitaire d'un groupe de Lie $G$. Pour toute distribution $\varphi\in \Cal E'(G)$ telle que~:
$$WF(\varphi)\cap -WF_\pi=\emptyset$$
l'op\'erateur $\pi(\varphi)$ est born\'e.
\dem
Ce corollaire r\'esulte du lemme suivant~:
\lemme{III.6}
Soit $A$ un op\'erateur d\'efini sur un domaine dense $\Cal D$ d'un espace de Hilbert $\Cal H$ tel que pour tout $T\in J_2$ laissant $\Cal D$ invariant l'op\'erateur $AT$ soit encore dans $J_2$. Alors $A$ est born\'e.
\dem
Supposons que $A$ ne soit pas born\'e. Il existe alors une suite de vecteurs $(\xi_j)$ de norme $1$ dans le domaine de $A$ telle que $A\xi_j=\eta_j$ avec $\|\eta_j\|=\alpha_j\rightarrow +\infty$. Quitte \`a prendre une sous-suite on peut supposer que l'on a~:
$$\alpha_j\ge j.$$
Soit $(e_j)$ une base orthonorm\'ee de $\Cal H$ contenue dans $\Cal D$, et soit l'op\'erateur $T:\Cal D\rightarrow \Cal D$ d\'efini par~:
$$Te_j=\beta_j\xi_j.$$
On suppose que $\sum \beta_j^2<+\infty$, de sorte que l'op\'erateur $T$ appartient \`a $J_2$. On a $ATe_j=\beta_j\eta_j$. Supposant $\beta_j=\frac 1j$, on a~:
$$\|AT\|_2^2=\sum_j\alpha_j^2|\beta_j|^2=+\infty,$$
autrement dit $AT$ n'appartient pas \`a $J_2$.
\qed
{\sl Remarque\/}~: ceci ne montre pas la continuit\'e de la correspondance $\varphi\mapsto \pi(\varphi)$ de $\Cal E'_\Gamma(G)$ dans $J_\infty(\pi)$. 
\th{III.7}
Soit $\pi$ une repr\'esentation unitaire d'un groupe de Lie $G$. Pour toute distribution $\varphi\in \Cal E'(G)$ telle que~:
$$WF(\varphi)\cap -WF_\pi=\emptyset$$
l'op\'erateur $\pi(\varphi)$ est r\'egularisant. De plus si $\pi$ est fortement tra\c cable, $\pi(\varphi)$ est un op\'erateur \`a trace. Pour toute partie conique ferm\'ee $\Gamma$ de $T^*G$ priv\'e de la section nulle et v\'erifiant~:
$$\Gamma\cap -WF_\pi=\emptyset,$$
la correspondance $\varphi\rightarrow \pi(\varphi)$ est s\'equentiellement continue de $\Cal E'_\Gamma(G)$ dans $J_1$. Enfin on a~:
$$\mop{Tr}\pi(\varphi)=<\chi_\pi.\varphi,\,1>.$$
\dem
Il suffit d'appliquer le corollaire III.5 \`a $v*\varphi$ pour tout $v\in\Cal U(\g g)$ pour voir que $\pi(\varphi)$ est bien r\'egularisant. Enfin si $\pi$ est fortement tra\c cable l'\'ecriture~:
$$A=\pi(1+\Delta)^{-s}\bigl(\pi(1+\Delta)^sA\bigr).$$  
montre que tous les op\'erateurs r\'egularisants sont \`a trace. On a de plus [DS lemma XI.9.20]~:
$$\|\pi(\varphi)\|_1\le \|\pi(\varphi*(1+\Delta)^{2s})\pi(1+\Delta)^{-s} \|_2 \|\pi(1+\Delta)^{-s} \|_2$$
En remarquant que $\varphi\mapsto \varphi*(1+\Delta)^{2s}$ est continue de $\Cal E'_\Gamma(G)$ dans $\Cal E'_\Gamma(G)$ et en appliquant le th\'eor\`eme III.4 avec $T=\pi(1+\Delta)^{-s}$ pour $s$ assez grand on obtient la continuit\'e s\'equentielle de $\varphi\mapsto \pi(\varphi)$. Enfin la formule \'enonc\'ee est imm\'ediate lorsque la distribution $\varphi$ appartient \`a $C^\infty_c(G)$, et $C^\infty_c(G)$ est s\'equentiellement dense dans $\Cal E'_\Gamma(G)$, ce qui permet d'\'etendre la formule \`a toutes les distributions $\varphi\in\Cal E'_\Gamma(G)$ si $\Gamma\cap -WF_\pi=\emptyset$.
\qed
\paragraphe{IV. Extension de la formule des caract\`eres}
\qquad
On suppose maintenant que la repr\'esentation $\pi$ est irr\'eductible, et qu'elle v\'erifie l'hypoth\`ese suivante~:
\smallskip
{\bf Hypoth\`ese (H)\/}~: 
\par
{\sl
$\pi$ est associ\'ee \`a une orbite coadjointe ferm\'ee $\Omega$ par la m\'ethode des orbites de Kirillov, de la mani\`ere suivante~: il existe un voisinage exponentiel $U$ de $0$ dans $\g g$, un entier strictement positif $\kappa(\pi)$ et une fonction $P_\Omega$ analytique partout non nulle sur $U$ et telle que $P_\Omega(0)=1$, tels que la formule des caract\`eres~:
$$\mop{Tr}\pi(T)=\kappa(\pi) \int_\Omega \Cal F\bigl(P_\Omega\inver .j_G.(T\circ\exp)\bigr)(\omega)\, d\beta_\Omega(\omega)\leqno{(*)}$$
soit v\'erifi\'ee pour toute fonction $T$ de classe $C^\infty$ \`a support compact inclus dans $\exp U$.\/}
\allowbreak 
\alinea{IV.1. Commentaires sur l'hypoth\`ese (H)}
Il existe de nombreux cas o\`u l'hypoth\`ese (H) est v\'erifi\'ee~: 
\ssq
1. Lorsque $G$ est compact connexe toute repr\'esentation unitaire irr\'eductible $\pi$ v\' erifie (H) avec $\kappa(\pi)=1$. C'est \'egalement vrai lorsque $G$ est nilpotent connexe, et dans ce cas $P_\Omega=1$ [Kir].
\ssq
2. Lorsque $G$ est r\'esoluble connexe et simplement connexe, \`a toute orbite coadjointe $\Omega$ telle que si $f\in \Omega$ l'indice du stabilisateur r\'eduit $\overline {G(f)}$ dans le stabilisateur $G(f)$ soit fini [Pk] on associe [Kh 1] une famille de repr\'esentations factorielles normales $\rho$ de type I. L'indice d\'efini ci-dessus est de la forme $n^2, n\in\N^*$ [Cha~1], et la repr\'esentation $\rho$ est de la forme $n\pi$ o\`u $\pi$ est unitaire irr\'eductible. Si de plus l'orbite $\Omega=G.f$ est ferm\'ee et temp\'er\'ee la repr\'esentation $\pi$ v\'erifie (H) avec $\kappa(\pi)=n$. Ceci g\'en\'eralise les r\'esultats de M. Duflo [BCD chap. IX] sur les caract\`eres des repr\'esentations associ\'ees \`a une orbite enti\`ere, c'est-\`a-dire dans le cas $n=1$.
\ssq
3. L'hypoth\`ese (H) est v\'erifi\'ee pour les s\'eries discr\`ete et principale g\'en\'eralis\'ee d'un groupe r\'eductif [Ros]. On peut prendre dans ce cas pour fonction $P_\Omega$ la fonction~:
$$P(x) =\Bigl(\det {\mop{sh ad}\frac x2\over \mop{ad}\frac x2}\Bigr)^{\frac 12}.$$
Cette fonction convient aussi dans le cas r\'esoluble lorsque l'orbite est de dimension maximale [Kh 1 \S\ 4.2.1].
\ssq
4. Enfin l'hypoth\`ese (H) est v\'erifi\'ee avec la m\^eme fonction $P$ que ci-dessus pour les repr\'esentations unitaires irr\'eductibles d'un groupe de Lie g\'en\'eral construites par M. Duflo [Du] associ\'ees \`a une orbite coadjointe ferm\'ee, temp\'er\'ee et de dimension maximale [Kh~3]. La formule obtenue co\"\i ncide avec la formule de Rossmann dans le cas r\'eductif (voir 3. ci-dessus). 
\ssq
Pour d\'ecrire l'entier $\kappa(\pi)$ il nous faut rappeler bri\`evement la construction de M. Duflo [Du]~: \`a un \'el\'ement $f\in \g g^*$ on associe un certain rev\^etement d'ordre deux $\wt {G(f)}$ du stabilisateur $G(f)$ et on d\'esigne par $\varepsilon$ l'\'el\'ement non trivial de $\wt {G(f)}$ qui se projette sur l'\'el\'ement neutre de $G(f)$. On d\'esigne par $X^{\smop{irr}}(f)$ l'ensemble des classes de repr\'esentations unitaires irr\'eductibles $\tau$ de $\wt {G(f)}$ dont la diff\'erentielle est multiple de la restriction de $-if$ \`a $G(f)$ (le signe moins provient de nos conventions sur la transform\'ee de Fourier) et telles que $\tau(\varepsilon)=-1$. Si $X^{\smop{irr}}(f)$ est non vide on dit que $f$ est admissible.
\ssq
On dit que $f$ est bien polarisable s'il existe en $f$ une polarisation r\'esoluble complexe v\'erifiant la condition de Pukanszky. La construction de M. Duflo consiste \`a associer \`a $f\in\g g^*$ admissible et bien polarisable et \`a $\tau\in X^{\smop{irr}}(f)$ une classe $T_{f,\tau}$ de repr\'esentations unitaires irr\'eductibles de $G$. Dans le cas o\`u l'orbite de $f$ est ferm\'ee, temp\'er\'ee et de dimension maximale M.S. Khalgui [Kh 3] montre que si de plus $\tau\in X^{\smop{irr}}(f)$ est de dimension finie l'hypoth\`ese (H) est alors v\'erifi\'ee pour $T_{f,\tau}$ avec~:
$$\kappa(T_{f,\tau})=\mop{dim}\tau.$$  
\ssq
Khalgui [Kh 2 \S\ 8] a par ailleurs trouv\'e un exemple de repr\'esentation tra\c cable qui ne v\'erifie pas (H)~: pour une certaine repr\'esentation unitaire irr\'eductible d'un groupe $G$ produit semi-direct d'un groupe semi-simple compact $K$ par son alg\`ebre de Lie $\g k$ il exhibe une orbite $\Omega$ temp\'er\'ee telle que la formule des caract\`eres est v\'erifi\'ee, mais avec une fonction $P_\Omega$ qui ne peut pas \^etre $C^\infty$ en $0$.
\alinea{IV.2. Extension de la formule des caract\`eres}
On se pose le probl\`eme suivant~: peut-on \'etendre pour une repr\'esentation v\'erifiant (H) la validit\'e de la formule des caract\`eres au-del\`a de $C^\infty_c(U)$? Autrement dit pour quelles distributions $\varphi$ \`a support compact l'op\'erateur $\pi(\varphi)$ est-il \`a trace, et \`a quelles conditions la trace de l'op\'erateur $\pi(\varphi)$ est-elle donn\'ee par
l'int\'egrale de la transform\'ee de Fourier de $P_\Omega\inver.\wt \varphi$ sur l'orbite?  
\ssq
Comme l'orbite $\Omega$ est suppos\'ee temp\'er\'ee, on voit que si la transform\'ee de Fourier de $P_\Omega\inver.\wt \varphi$ d\'ecro\^\i t suffisamment rapidement sur l'orbite le membre de droite dans la formule des caract\`eres est convergent. Nous allons donc montrer que l'op\'erateur $\pi(\varphi)$ est \`a trace si $WF(\varphi)$ ne rencontre pas l'oppos\'e du c\^one asymptote \`a $\Omega$, et que dans ce cas $\mop{Tr}\pi(\varphi)$ est donn\'ee par la formule des caract\`eres (*).
\prop{IV.1}
Toute repr\'esentation $\pi$ unitaire irr\'eductible v\'erifiant l'hypoth\`ese (H) est fortement tra\c cable.
\dem
Le point cl\'e est le fait que l'orbite $\Omega$ associ\'ee est temp\'er\'ee [Cha~2,3]. On aura \'egalement besoin d'un lemme d'analyse fonctionnelle~:
\lemme{IV.2}
Soit $(A_k)_{k\ge 1}$ une suite d'op\'erateurs de Hilbert-Schmidt dans un espace de Hilbert s\'eparable $\Cal H$, et soit $A$ un op\'erateur born\'e, tels que~:
\smallskip
1) $A_k$ converge fortement vers $A$ lorsque $k\rightarrow +\infty$ 
\smallskip
2) $\|A_k\|_2\fleche 8_{k\rightarrow +\infty}\ell\in[0,+\infty[$.
\smallskip
Alors $A$ est un op\'erateur de Hilbert-Schmidt.
\dem
On se donne une base orthonorm\'ee $(e_i)$ de $\Cal H$. La suite de fonctions positives~:
$$\varphi_k:i\longmapsto \|A_ke_i\|^2$$
converge ponctuellement vers $\varphi:i\mapsto \|Ae_i\|^2$. Le lemme de Fatou permet de montrer que la s\'erie~:
$$\sum_i\|Ae_i\|^2$$
est convergente. D'apr\`es [DS lemma XI.9.32] l'op\'erateur $A$ est donc dans $J_2$.
\qed
{\sl Fin de la d\'emonstration de la proposition IV.1\/}~: soit $V$ un voisinage exponentiel de $e$ dans $G$. Soit $\varphi\in C^m_c(V)$ avec $m$ entier positif. Soit $(\alpha_k)$ une approximation de l'identit\'e dans $C^\infty_c(G)$, telle que pour $k$ assez grand le support de $\varphi_k=\alpha_k*\varphi$ soit encore contenu dans $V$. Pour tout entier pair $2j\le m-n-1$ la fonction $(1+\|\xi\|^2)^{j}\Cal F(P_\Omega\inver\wt\varphi_k)$ converge uniform\'ement sur $\g g^*$ vers $(1+\|\xi\|^2)^{j}\Cal F(P_\Omega\inver\wt\varphi)$. Comme l'orbite est temp\'er\'ee [Cha~2,3], si $m$ est assez grand il existe donc $j$ tel que l'int\'egrale~:
$$\int_\Omega \Cal F(P_\Omega\inver \wt\varphi)(\xi)d\beta_\Omega(\xi)
=\int_\Omega (1+\|\xi\|^2)^j\Cal F(P_\Omega\inver \wt\varphi)(\xi)\,(1+\|\xi\|^2)^{-j}d\beta_\Omega(\xi)$$
est convergente et la mesure $(1+\|\xi\|^2)^{-j} d\beta_\Omega$ est finie. En appliquant la formule des caract\`eres \`a $\pi(\varphi_k)$ et le th\'eor\`eme de convergence domin\'ee de Lebesgue on a~:
$$\mop{Tr}(\varphi_k)\fleche 8_{k\rightarrow +\infty}
\int_\Omega \Cal F(P_\Omega\inver \varphi)(\xi)d\beta_\Omega(\xi)$$
L'op\'erateur $\pi(\varphi)$ est born\'e puisque $\varphi\in L^1(G)$. Appliquant alors le lemme IV.2 \`a la suite d'op\'erateurs~:
$$A_k=\pi(\varphi_k),\ A=\pi(\varphi)$$
on voit qu'il existe un entier $m$ tel que pour tout $\varphi\in C^m_c(V)$ l'op\'erateur $\pi(\varphi)$ appartient \`a $J_2$. Quitte \`a augmenter sensiblement l'entier $m$ on peut alors passer de $J_2$ \`a $J_1$~: on utilise le r\'esultat suivant (M. Duflo dans [BCD lemme IX.3.2.3], cf. aussi [J\o] theorem 4), cons\'equence de l'ellipticit\'e au voisinage de l'\'el\'ement neutre de l'op\'erateur invariant \`a droite donn\'e par $\Delta$~: quitte \`a r\'etr\'ecir un peu l'ouvert $V$ il existe un entier $m'$, une distribution $\beta\in \Cal E'(V)$ et une fonction $\gamma\in C^\infty_c(V)$ tels que~:
$$\Delta^{*m'}*\beta=\delta+\gamma$$
o\`u $\delta$ d\'esigne la mesure de Dirac en l'\'el\'ement neutre. L'op\'erateur invariant \`a droite donn\'e par $\Delta^{m'}*_-$ \'etant elliptique d'ordre $2m'$, la distribution $\beta$ est une fonction de classe $C^{2m'-n-1}$, donc de classe $C^m$ si $m'$ est assez grand. On d\'eduit de ceci, comme dans [BCD p. 251-252] que l'on a~:
$$\varphi=\varphi*\Delta^{*2m'}*\beta*\beta-\varphi*\Delta^{*m'}*(\gamma *\beta +\beta *\gamma)
	+\varphi*\gamma*\gamma$$
d'o\`u on d\'eduit~:
$$\pi(\varphi)=\pi(\varphi)\circ\Bigl(\pi(\Delta)^{2m'}\circ\pi(\beta)^2
+\pi(\Delta)^{m'}\circ\bigl(\pi(\gamma)\pi(\beta)+\pi(\beta)\pi(\gamma)\bigr)
+ \pi(\gamma)^2\Bigr)$$
Si $\varphi\in C^{m+4m'}_c(G)$ les op\'erateurs $\pi(\varphi)\circ\pi(\Delta)^{2m'}$, $\pi(\varphi)\circ\pi(\Delta)^{m'}$ et $\pi(\varphi)$ sont dans $J_2$ d'apr\`es ce qui pr\'ec\`ede, donc born\'es, et les op\'erateurs $\pi(\beta)^2$, $\pi(\gamma)\pi(\beta)+\pi(\beta)\pi(\gamma)$ et $\pi(\gamma)^2$ sont \`a trace puisque $\pi(\beta)$ et $\pi(\gamma)$ sont dans $J_2$. Donc $\pi(\varphi)$ est \`a trace, ce qui d\'emontre la proposition, d'apr\`es le crit\`ere 3) de la proposition I.1. 
\qed
\th{IV.3}
Soit $\pi$ une repr\'esentation unitaire irr\'eductible de $G$ v\'erifiant l'hypoth\`ese (H). On d\'esigne par $AC(\Omega)$ le c\^one asymptote \`a l'orbite, c'est-\`a-dire l'ensemble des $\xi\in\g g^*$ tels que tout voisinage conique de $\xi$ a une intersection non born\'ee avec $\Omega$. Soit $V$ un voisinage exponentiel de $e$ dans $G$. Alors~:
\smallskip
1) $WF_\pi^0\subset -AC(\Omega)$
\smallskip
2) Pour toute distribution $\varphi\in \Cal E'(V)$ telle que~:
$$WF(\varphi)\cap -WF_\pi=\emptyset,$$
l'op\'erateur $\pi(\varphi)$ est \`a trace, et la formule des caract\`eres~:
$$\mop{Tr}\pi(\varphi)=\kappa(\pi) \int_\Omega \Cal F(P_\Omega\inver .\wt \varphi)(\omega)\, d\beta_\Omega(\omega)\leqno{(*)}$$
est v\'erifi\'ee.
\dem
Le 1) est d\^u \`a R. Howe [Hw proposition 2.2], qui en donne une d\'emonstration succincte dans le cas nilpotent. Pour une d\'emonstration dans le cadre g\'en\'eral de la transform\'ee de Fourier d'une mesure positive temp\'er\'ee quelconque, voir [D-V lemme 4]. Rappelons-en le principe~: au vu de la proposition II.2 et de la formule des caract\`eres un point $\xi_0\in\g g^*$ n'appartient pas \`a $WF_\pi^0=WF^0(\chi_\pi)$ si et seulement si il existe un voisinage $W$ de $\xi_0$ et un voisinage ouvert $U$ de $e$ dans $G$ tels que pour tout $\xi\in W$ et toute fonction $\alpha\in C^\infty_c(U)$ l'int\'egrale~:
$$\int_\Omega\Cal F(P_\Omega\inver \wt\alpha)(\omega+t\xi)d\beta_\Omega\leqno{(I)}$$
est \`a d\'ecroissance rapide lorsque $t\rightarrow +\infty$, et ce de mani\`ere uniforme lorsque $\xi$ varie dans un voisinage compact de $\xi_0$ dans $W$. Supposons que $\xi_0$ n'appartienne pas \`a $-AC(\Omega)$. Alors il existe un voisinage $W$ de $\xi_0$ et une constante $B>0$ telle que pour tout $\omega\in \Omega$ et pour tout $t$ assez grand on ait~:
$$\|\omega+t\xi\|\ge Bt\|\xi\|.$$
Comme par ailleurs la transform\'ee de Fourier $\Cal F(P_\Omega\inver \wt\alpha)$ est \`a d\'ecroissance rapide on montre facilement en utilisant le fait que la mesure $d\beta_\Omega$ est temp\'er\'ee que $(I)$ est \`a d\'ecroissance rapide en $t$.
\ssq
Pour le 2) il suffit de remarquer que le membre de droite de la formule (*) se prolonge en une forme lin\'eaire continue sur $\Cal E'_\Gamma(V)$ pour toute partie conique ferm\'ee $\Gamma$ de $T^*V\backslash 0$ ne rencontrant pas $WF_\pi$. Par ailleurs d'apr\`es le th\'eor\`eme III.7 le membre de gauche d\'efinit une forme lin\'eaire s\'equentiellement continue sur $\Cal E'_\Gamma(G)$, et donc aussi sur $\Cal E'_\Gamma(V)$. Ces deux formes lin\'eaires sont donc \'egales puisqu'elles sont \'egales sur le sous-espace s\'equentiellement dense $C^\infty_c(V)$ ([Hr 1\S\ 2.5]).
\qed
{\sl Remarque\/}~: la proposition 2.2 de [Hw] stipule l'\'egalit\'e de $WF_\pi^0$ et de $-AC(\Omega)$ dans le cas nilpotent. L'id\'ee de d\'emonstration est la suivante : on suppose $G$ connexe et simplement connexe, on identifie $G$ et son alg\`ebre de Lie via l'exponentielle et on consid\`ere la famille de fonctions ``presque \`a support compact'' $(\varphi_\tau)_{\tau\in [2 ,+\infty]}$ d\'efinie par~:
$$\varphi_\tau(x)=\tau\inver e^{-{\|(\log \tau)x\|^2\over 2}}$$
pour $\tau<+\infty$ et $\varphi_\infty=0$. Cette famille de fonctions forme un compact de $C^\infty(\g g)$. On montre alors que pour $\xi\in -AC(\Omega)$ il existe une suite $\xi_k\rightarrow \xi$ telle que l'expression~:
$$\Cal F(\varphi_\tau.\chi_\pi)(\tau \xi_k)$$
n'est pas \`a d\'ecroissance rapide en $\tau$. Il suffit de consid\'erer une suite $\xi_k\rightarrow \xi$ telle qu'il existe $\tau_k\rightarrow +\infty$ tel que $-\tau_k\xi_k\in \Omega$ (par d\'efinition du c\^one asymptote). Il existe un \'el\'ement $\omega_0\in \Omega$ et une suite $g_k\in G$ telle que $-\tau_k\xi_k=\mop{Ad}^*\omega_0$. On a alors~:
$$\eqalign{I_k:=\Cal F(\varphi_{\tau_k}.\chi_\pi)(\tau_k \xi_k)	&= \int_\Omega
	\Cal F.\varphi_{\tau_k}(\omega+\tau_k\xi_k)\,d\beta_\Omega(\omega)\cr
							&= \int_\Omega
	\Cal F.\varphi_{\tau_k}\bigl(\mop{Ad}^*g_k.(\omega-\omega_0)\bigr)\,d\beta_\Omega(\omega).\cr}$$
Il existe un voisinage compact $V$ de $0$ dans $\g g^*$ tel que~:
$$\int_{\{\omega\in\Omega/\omega-\omega_0\in V \}}d\beta_\Omega(\omega)\ge 2$$
et un $\tau_0$ tel que tel que~:
$$\mopl{inf}_{\eta\in V}(\Cal F.\varphi_\tau)(\eta)\ge{\tau\inver(\log \tau)^{-\frac n2}\over 2}.$$
On a donc~:
$$\eqalign{I_k	&\ge{\tau_k\inver(\log \tau_k)^{-\frac n2}\over 2} \int_{\omega-\omega_0\in\smop{Ad}^*g_k\inver.V }d\beta_\Omega(\omega)\cr
		&\ge {\tau_k\inver(\log \tau_k)^{-\frac n2}\over 2}\|\mop{Ad}^*g_k\|^{-\smop{dim}\Omega} \int_{\omega-\omega_0\in V }d\beta_\Omega(\omega)\cr
		&\ge \tau_k\inver(\log \tau_k)^{-\frac n2}\|\mop{Ad}^* g_k\|^{-\smop{dim}\Omega}.\cr}$$
Par nilpotence de $G$ on peut choisir $g_k$ tel que $\|\mop{Ad}^*g_k\|$ soit \`a croissance polynomiale en $\tau_k$, ce qui montre que $I_k$ n'est pas \`a d\'ecroissance rapide en $\tau_k$. 
On se donne alors une fonction $\gamma\in C^\infty_c(U)$ \'egale \`a $1$ au voisinage de $0$ dans $\g g$, o\`u $U$ est un voisinage quelconque de $0$. Si on pose~:
$$\varphi'_\tau=\gamma\varphi_\tau$$
on v\'erifie facilement que~:
$$\mopl{sup}_x(\varphi_\tau(x)-\varphi'_\tau(x))=O(\tau^{-N})$$
pour tout $N$. L'int\'egrale~:
$$I'_k:=\Cal F(\varphi'_{\tau_k}.\chi_\pi)(\tau_k \xi_k)$$
n'est donc pas \`a d\'ecroissance rapide en $\tau$. Le vecteur $\xi$ appartient donc au front d'onde de $\chi_\pi$, donc \`a $WF^0_\pi$ d'apr\`es la proposition II.2.
\paragraphe{V. Applications}
\alinea{V.1. Application aux op\'erateurs pseudo-diff\'erentiels}
On introduit dans [M1, M2] les classes de symboles analytiques de la fa\c con suivante~: pour tout voisinage compact $Q$ de $0$ dans $\g g$, pour tout r\'eel $m$ et pour tout $\rho\in [0,1]$ on d\'efinit la classe $AS_\rho^{m,Q}(\g g^*)$ comme l'espace des fonctions $p$ sur $\g g^*$ dont le support de la transform\'ee de Fourier inverse est inclus dans $Q$, et qui v\'erifient les estimations~:
$$|D^\alpha p(\xi)|\le C_\alpha(1+\|\xi\|)^{m-\rho|\alpha|}.$$
Les symboles utilis\'es sont des transform\'ees de Fourier de distributions \`a support compact sur $\g g$, dont le support singulier est r\'eduit \`a $\{0\}$ lorsque $\rho>0$. Le cadre des distributions \`a support compact que nous avons choisi ici correspond par transform\'ee de Fourier \`a la classe $AS^{M,Q}_\rho$ avec $\rho=0$. La convolution sur le groupe donne donc naissance \`a un produit~:
$$\#:AS^{M,Q}_0\times AS^{M,Q}_0\longrightarrow AS^{M,Q^2}_0$$
si $Q$ est assez petit. La restriction de ce produit \`a $AS^{M,Q}_\rho\times AS^{M,Q}_\rho$ pour $\rho>1/2$ permet de d\'evelopper [M1] un calcul symbolique similaire \`a plusieurs \'egards \`a celui des op\'erateurs pseudo-diff\'erentiels sur $\R^n$. Pour toute repr\'esentation unitaire $\pi$ de $G$ on d\'efinit alors l'op\'erateur pseudo-diff\'erentiel de symbole $p$ dans l'espace de la repr\'esentation $\pi$ par~:
$$p^{W,\pi}=\pi(\varphi)$$
avec $j_G(x)\varphi(\exp x)=\wt \varphi(x)=\Cal F\inver p(x)$ pour $x\in\g g$. Le support singulier de $\varphi$ \'etant r\'eduit \`a l'\'el\'ement neutre son front d'onde s'identifie \`a un c\^one dans $\g g^*$, et l'application du th\'eor\`eme III.7 et du th\'eor\`eme IV.3.1) dans ce cas particulier nous redonne le th\'eor\`eme III.1 de [M2].
\alinea{V.2. Restriction \`a un sous-groupe ferm\'e transverse}
\qquad
On se donne un groupe de Lie $G$ d'alg\`ebre de Lie $\g g$, et une 
repr\'esentation unitaire fortement tra\c cable $\pi$ de $G$. On se donne
un sous-groupe ferm\'e $H$ transverse au front d'onde de $\pi$, c'est \`a dire
tel que $WF^0_\pi\cap \g h^\perp=\emptyset$. On se fixe une mesure de Haar \`a gauche $dh$ sur $H$
\ssq
Toute fonction $\psi\in C^\infty_c(H)$ donne naissance \`a une distribution $\overline \psi$ \`a support compact sur $G$, d\'efinie par~:
$$\eqalign {<\overline \psi,\, \varphi>	&=<\psi,\, \varphi\restr H>\cr
	&=\int_{H}\psi(h)\varphi(h)\,dh.	\cr}$$
On a bien entendu~:
$$WF_e(\overline \psi)\subset \g h^\perp.$$
On en d\'eduit facilement que $WF(\overline\psi)\cap -WF_\pi=\emptyset$.
On utilise alors le th\'eor\`eme III.7 : l'op\'erateur $\pi(\overline\psi)$ est \`a trace. Or on a l'\'egalit\'e entre op\'erateurs~:
$$\pi(\overline\psi)=\pi\restr H (\psi).$$
La repr\'esentation $\pi\restr H$ est donc tra\c cable, et on a bien~:

$$\eqalign{<\chi_{(\pi\srestr H)},\,\psi>
	&=Tr(\pi\restr H(\psi))	\cr
	&=Tr(\pi(\overline \psi))	\cr
	&=<\chi_\pi.\overline \psi,\, 1>
		\hbox{ d'apr\`es le th\'eor\`eme III.7}	\cr
	&=<(\chi_\pi)\restr H,\, \psi>
		\hbox{ par d\'efinition de la restriction d'une distribution.}\cr}$$
d'o\`u l'\'egalit\'e~:
$$\chi_{(\pi\srestr H)}=(\chi_\pi)\restr H.$$
Autrement dit~:
\th{V.1}
Soit $G$ un groupe de Lie d'alg\`ebre de Lie $\g g$, soit $\pi$ une repr\'esentation unitaire fortement tra\c cable de $G$, et soit $H$ un sous-groupe ferm\'e de $G$ d'alg\`ebre de Lie $\g h$ v\'erifiant la condition~:
$$WF^0_\pi\cap \g h^\perp=\emptyset.$$
Alors la restriction de $\pi$ au sous-groupe $H$ est tra\c cable, et le caract\`ere de $\pi\restr H$ est donn\'e par la restriction du caract\`ere de $\pi$ au sous-groupe $H$.
\ndem
\cor{V.2}
Sous les hypoth\`eses du th\'eor\`eme V.1, la restriction de $\pi$ au sous-groupe $H$ est somme directe de repr\'esentations unitaires irr\'eductibles tra\c cables de $H$, chacune intervenant avec multiplicit\'e finie.
\ndem
Ce r\'esultat a \'et\'e conjectur\'e par M. Duflo et M. Vergne ([D-V], fin du chapitre 1).

\bigskip
\paragraphe{R\'ef\'erences :}
\bib{BCD}P. Bernat, N. Conze, M. Duflo et al. {\sl Repr\'esentations des groupes de Lie r\'esolubles\/}, Dunod 1972.
\bib{Cha 1}J-Y. Charbonnel, {\sl La mesure de Plancherel pour les groupes de lie r\'esolubles connexes\/}, Lect. Notes in Math. 576 (1977), 32-76.
\bib{Cha 2}J-Y. Charbonnel, {\sl Orbites ferm\'ees et orbites temp\'er\'ees\/}, Ann. Sci. Ec. Norm. Sup. 23 (1990), 123-149.
\bib{Cha 3}J-Y. Charbonnel, {\sl Orbites ferm\'ees et orbites temp\'er\'ees, II\/}, J. Funct. Anal. 138 (1996), 213-222. 
\bib{Du}M. Duflo, {\sl Construction de repr\'esentations unitaires d'un groupe de Lie\/}, Cours d'\'et\'e du C.I.M.E, Cortona 1980, 129-221.
\bib{D-V}M. Duflo, M. Vergne, {\sl Orbites coadjointes et cohomologie \'equivariante\/}, in {\sl The orbit method in representation theory\/}, Copenhague 1988. Progr. Math. 82, M. Duflo, N.V. Pedersen, M. Vergne eds.
\bib{D-S}N. Dunford, J.T. Schwartz, {\sl linear operators\/}, vol II, Interscience, New-York 1958.
\bib{Dui}J.J. Duistermaat, {\sl Fourier integral operators\/}, Courant Institute, 1973 (R\'e\'ed. Progress in Math. 130, Birkh\"auser 1995).
\bib{Ga}A. Gabor, {\sl Remarks on the wave front set of a distribution\/}, Trans. Amer. Math. Soc. 170 (1972), 239-244.
\bib{Go}R. Goodman, {\sl Elliptic and subelliptic estimates for operators in an enveloping algebra\/}, Duke Math. J. 47, No 4 (1980), 819-833.
\bib{Hei}D.B. Heifetz, {\sl p-adic oscillatory integrals and wave front sets\/}, Pac. J. Math. 116 No 2 (1985), 285-305.
\bib{Hr 1}L. H\"ormander, {\sl Fourier integral operators I\/}, Acta Math. 127 (1971), 79-183.
\bib{Hr 2}L. H\"ormander, {\sl The Weyl calculus of pseudodifferential operators,\/} Comm. Pure Appl. Math. 32 (1979), 359-443.
\bib{Hr 3}L. H\"ormander, {\sl The analysis of partial differential operators I\/}, Springer 1983.
\bib{Hu}A. Hulanicki, {\sl Subalgebra of $L_1(G)$ associated with laplacian on a Lie group\/}, Coll. Math. 31 No 2 (1974), 259-287.
\bib{Hw}R. Howe, {\sl Wave front sets of representations of Lie groups\/}, in {\sl Automorphic forms, representation theory and Arithmetic\/}, Tata Inst. fund. res. st. math. 10, Bombay 1981.
\bib{J\o}P.E.T. J\o rgensen, {\sl Distribution representations of Lie groups\/}, J. Math. Anal. Appl. 65 (1978), 1-19.
\bib{Kh 1}M.S. Khalgui, {\sl Sur les caract\`eres des groupes de Lie r\'esolubles\/}, Publ. Math. Univ. Paris VII 2 (1978).
\bib{Kh 2}M.S. Khalgui, {\sl Sur les caract\`eres des groupes de Lie \`a radical cocompact\/}, Bull. Soc. Math. France 109 (1981), 331-372.
\bib{Kh 3}M.S. Khalgui, {\sl Caract\`eres des groupes de Lie\/}, J. Funct. Anal. 47 (1982), 64-77 
\bib{Kir}A.A. Kirillov, {\sl Elements of the theory of representations\/}, Springer 1976.
\bib{M1}D. Manchon, {\sl Weyl symbolic calculus on any Lie group\/}, Acta Appl. Math. 30 (1993), 159-186.
\bib{M2}D. Manchon, {\sl op\'erateurs pseudodiff\'erentiels et repr\'esentations unitaires des groupes de Lie\/}, Bull. Soc. Math. France 123 (1995), 117-138.
\bib{Ne}E. Nelson, {\sl Analytic vectors\/}, Ann. Math. 70 No 3 (1959), 572-615.
\bib{Ne-St}E. Nelson, W.F. Stinespring, {\sl Representation of elliptic operators in an enveloping algebra\/}, Amer. J. Math. 81 (1959), 547-560.
\bib{Pk}L. Pukanszky, {\sl Unitary representations of solvable Lie groups\/}, Ann. Ec. Norm. Sup. 4 (1971), 457-608.
\bib{Rob}D. W. Robinson, {\sl Elliptic operators and Lie groups\/}, Oxford 1991. 
\bib{Ros}W. Rossmann, {\sl Kirillov's character formula for reductive Lie groups\/}, Invent. Math. 48 (1978), 207-220.
\bib{Shu}M.A. Shubin, {\sl Pseudodifferential operators and spectral theory\/}, Springer 1987.
\bye